\documentclass[a4paper,twoside,12pt]{article}
\usepackage{amssymb} 
\begin{document}
\title{The extremal genus embedding of  graphs \footnotemark[2]
\author{Guanghua Dong$^{1,2}$, Han Ren$^{3}$, Ning Wang$^{4}$, Hao Wu$^{3}$\\
{\small\em 1.Department of Mathematics, Normal
University of Hunan, Changsha, 410081, China}\\
\hspace{-1mm}{\small\em 2.Department of Mathematics, Tianjin
Polytechnic
University, Tianjin, 300387, China}\\
\hspace{-1mm}{\small\em 3.Department of Mathematics, East China
Normal University, Shanghai, 200062,China}\\
\hspace{-5mm}{\small\em 4.Department of Information Science and
Technology,
Tianjin University of Finance }\\
\hspace{-74mm} {\small\em and Economics, Tianjin, 300222, China}\\
}} \footnotetext[2]{\footnotesize \em   This work was partially
Supported  by the China Postdoctoral Science Foundation funded
project (Grant No: 20110491248), the National Natural Science
Foundation of China (Grant No: 11171114), and  the New Century
Excellent Talents in University (Grant No: NCET-07-0276).}
 \footnotetext[1]{\footnotesize \em E-mail: gh.dong@163.com(G. Dong); hren@math.ecnu.edu.cn(H. Ren);  ninglw@163.com(N. Wang). }

\date{}
\maketitle

\vspace{-.8cm}
\begin{abstract}

Let $W_{n}$ be a wheel graph with $n$ spokes. How does the genus
change if adding a degree-3 vertex $v$, which is not in $V(W_{n})$,
to the graph $W_{n}$? In this paper, through the joint-tree model we
obtain that the genus of $W_{n}+v$ equals 0 if the three neighbors
of $v$ are in the same face boundary  of $\mathbb{P}(W_{n})$;
otherwise, $\gamma(W_{n}+v)=1$, where $\mathbb{P}(W_{n})$ is the
unique planar embedding of $W_{n}$. In addition, via the independent
set, we provide a lower bound on the maximum genus of graphs, which
may be better than both the result of D. Li \& Y. Liu and the result
of Z. Ouyang $etc.$ in Europ. J. Combinatorics. Furthermore, we
obtain a relation between the independence number and the maximum
genus of graphs, and provide an algorithm to obtain the lower bound
on the number of the distinct maximum genus embedding of the
complete graph $K_{m}$, which, in some sense, improves the result of
Y. Caro and S. Stahl respectively.

\bigskip
\noindent{\bf Key Words:} joint-tree model;  genus; maximum genus; independence number \\
{\bf MSC(2000):} \ 05C10
\end{abstract}


\bigskip
\noindent {\bf 1. Introduction}

Graph considered here are all finite and connected. If the graph $M$
can be obtained from a graph $G$ by successively contracting edges
and deleting edges and isolated vertices, then $M$ is a $minor$ of
$G$. The minimum genus $\gamma_{min}(G)$ (or, simply, the genus
$\gamma(G)$) of a graph $G$ is the minimum integer $g$ such that
there exists an embedding of $G$ into the orientable surface $S_{g}$
of genus $g$, and the \emph{maximum genus} $\gamma_M(G)$ of a
connected graph \emph{G} is the maximum integer \emph{k} such that
there exists an embedding of $G$ into the orientable surface of
genus $k$. The difference between  the maximum genus and the minimum
genus of a graph $G$ is called the $genus$ $range$ of $G$. A graph
$G$ is said to be \emph{upper embeddable} if $\gamma_M(\emph{G})$ =
$\lfloor\frac{\beta(G)}{2}\rfloor$, where $\beta(G)$ is the $cycle$
$rank$ (or $Betti$ $number$) of $G$. A $one$-$face$ embedding
($two$-$face$ embedding) $\psi(G)$ of a graph $G$ means that the
face number of $\psi(G)$ is one (two). An $odd$ $vertex$ is a vertex
whose degree is an odd number. For $n\geqslant 3$, the $wheel$ of
$n$ spokes is the graph $W_{n}$ obtained from the $n$-cycle $C_{n}$
by adding a new vertex (called the $center$ of the $wheel$) and
joint it to all vertices of $C_{n}$. For example, $W_{3}=K_4$. A
$subdivision$ of an edge $e\in E(W_{n})$ means inserting a vertex of
degree two to $e$, where the inserted vertex is called a
$subdividing$-$vertex$ of $W_{n}$. Let $v$ be a degree-three vertex
which is not in $V(W_{n})$, then the graph $W_{n}+v$, which is
called the $near$-$wheel$ graph, means the connected graph obtained
from $W_{n}$ by joining $v$ to $v_{i} \ (i=1, 2, 3)$, where $v_{i}$
may be a $subdividing$-$vertex$ of $W_{n}$ or a vertex which belongs
to $V(W_{n})$. Furthermore, the
 vertices $v_1$, $v_2$, $v_3$ are called the
$antennal$-$vertex$ of the graph $W_{n}+v$.

Surfaces considered here are compact 2-dimensional manifold without
boundary. An orientable surface $S$ can be regarded as a polygon
with even number of directed edges such that both $a$ and $a^{-}$
occurs once on $S$ for each $a\in S$, where the power
\textquotedblleft $-$\textquotedblright  means that the direction of
$a^{-}$ is opposite to that of $a$ on the polygon. For convenience,
a polygon is represented by a linear sequence of lowercase letters.
An elementary result in algebraic topology states that each
orientable surface is equivalent to one of the following standard
forms of surfaces:
\[ O_{p}=\left\{
 \begin{array}{ll}
 a_{0}a_{0}^{-},  &\mbox{$p = 0$,}\\
 \prod\limits_{i=1}^{p}a_{i}b_{i}a_{i}^{-}b_{i}^{-},  &\mbox{$p \geq 1$ .}
 \end{array}
 \right.
\]
which are the sphere ($p=0$), torus ($p=1$), and the orientable
surfaces of genus $p\ (p\geq2)$. The genus of a surface $S$ is
denoted by $g(S)$. Let  $A$, $B$, $C$, $D$, and $E$ be possibly
empty linear sequence of letters. Suppose $A=a_1a_2\dots a_{r},
r\geq1$, then $A^{-}=a_{r}^{-}\dots a_2^{-}a_1^{-}$ is called the
$inverse$ of $A$. If  $\{a, b, a^{-}, b^{-}\}$ appear in a sequence
of the form of $AaBbCa^{-}Db^{-}E$, then they are said to be an
$interlaced$ $set$; otherwise, a $parallel$ $set$. Let
$\widetilde{S}$ be the set of all surfaces. For a surface $S\in
\widetilde{S}$, we obtain its genus $g(S)$ by using the following
transforms to determine its equivalence to one of the standard
forms.

\medskip

{\bf Transform 1} \ \ $Aaa^{-} \sim A$, where $A\in \widetilde{S}$
and $a\notin A$.

{\bf Transform 2} \ \  $AabBb^{-}a^{-} \sim AcBc^{-}$.

{\bf Transform 3} \  \ $(Aa)(a^{-}B) \sim (AB)$.

{\bf Transform 4} \ \  $AaBbCa^{-}Db^{-}E \sim ADCBEaba^{-}b^{-}$.

\medskip

\noindent   In the above transforms, the parentheses stand for
cyclic order. For convenience, the parentheses are always omitted
when unnecessary to distinguish cyclic or linear order. For more
details concerning surfaces, the reader is referred to \cite{liu1}
and \cite{rin}.

Let $T$ be a spanning tree of a graph $G=(V,E)$, then
$E=E_{T}+\hat{E}_{T}$, where $E_{T}$ consists of all the tree edges,
and $\hat{E}_{T}=\{ \hat{e}_1, \hat{e}_2, \dots \hat{e}_{\beta}\}$
consists of all the co-tree edges, where $\beta=\beta(G)$ is the
cycle rank of $G$. Split each co-tree edge
$\hat{e}_{i}=(u[\hat{e}_{i}], v[\hat{e}_{i}]) \in \hat{E}_{T}$ into
 two semi-edges $(u[\hat{e}_{i}], v_{i})$,
$(v[\hat{e}_{i}], \bar{v}_{i})$, denoted by $\hat{e}_{i}^{+}$ and
$\hat{e}_{i}^{-}$ respectively.  Let $\widetilde{T}=(V+V_1, E+E_1)$,
where $V_1=\{v_{i}, \bar{v}_{i} | 1\leq i \leq \beta \}$,
$E_1=\{(u[\hat{e}_{i}], v_{i}), (v[\hat{e}_{i}], \bar{v}_{i}) |
1\leq i \leq \beta\}$. Obviously, $\widetilde{T}$ is a tree. A
rotation at a vertex $v$, which is denoted by $\sigma_{v}$, is a
cyclic permutation of edges incident on $v$. A  rotation system
$\sigma=\sigma_{G}$ for a graph $G$  is a set $\{\sigma_{v} |
\forall v \in V(G)\}$. The tree $\widetilde{T}$ with a rotation
system of $G$ is called a $joint$-$tree$ of $G$, and is denoted by
$\widetilde{T}_{\sigma}$. Because $\widetilde{T}_{\sigma}$ is a
tree, it can be embedded in the plane. By reading the lettered
semi-edges of $\widetilde{T}_{\sigma}$ in a fixed direction
(clockwise or anticlockwise), we can get an algebraic representation
of the surface which is represented by a $2\beta-$polygon. Such a
surface, which is denoted by $S_{\sigma}$, is called an associated
surface of $\widetilde{T}_{\sigma}$. A joint-tree
$\widetilde{T}_{\sigma}$ of $G$ and its associated surface is
illustrated by Fig.1, where the rotation at each vertex of $G$
complies with the clockwise rotation. From \cite{liu1}, there is 1-1
correspondence between the associated surfaces (or joint-trees) and
the embeddings of a graph. The $joint$-$tree$  is originated from
the early works of Liu \cite{liu2}, and more detailed information
about the $joint$-$tree$ can be found in \cite{liu1}. Terminologies
and notations not defined here can be seen in \cite{bon} for graph
theory and \cite{moh} for topological graph theory.


\bigskip
\setlength{\unitlength}{1mm}
\begin{center}
\begin{picture}(100,30)

\put(15,3){\begin{picture}(10,10)

\put(29,-5){{\bf Fig. 1.}}
\end{picture}
}


\put(-15,4)

{\begin{picture}(10,10)

\put(0,10){\circle*{1.5}}

\put(26,10){\circle*{1.5}}

\put(13,28){\circle*{1.5}}

\put(13,17){\circle*{1.5}}

\put(0.8,10.5){\line(2,1){11.6}}

\put(25,10.5){\line(-2,1){11.6}}

\put(13,27){\line(0,-1){9.5}}

\put(0.7,11){\line(3,4){12}}

\put(25.7,11){\line(-3,4){12}}

\put(1,10){\line(1,0){24}}

\put(0.5,20){ $\hat{e}_1$}

\put(20,20){ $\hat{e}_2$}

\put(9,6){ $\hat{e}_3$}

\put(12,0){ $G$}

\end{picture}

}

\put(35,4){\begin{picture}(10,10)

\put(7,14){\circle*{1.5}}

\put(19,14){\circle*{1.5}}

\put(13,23.2){\circle*{1.5}}

\put(13,17){\circle*{1.5}}

\put(7,14){\line(2,1){6}}

\put(19,14){\line(-2,1){6}}

\put(13,23.2){\line(0,-1){6}}

\put(7,14){\line(-1,1){4}}

\put(7,14){\line(-1,-1){4}}

\put(19,14){\line(1,1){4}}

\put(19,14){\line(1,-1){4}}

\put(13,23.2){\line(-1,1){4}}

\put(13,23.2){\line(1,1){4}}

\put(-3,18){ $\hat{e}_1$}

\put(3,28){ $\hat{e}_1^{-}$}

\put(17,28){ $\hat{e}_2$}

\put(24,18){ $\hat{e}_2^{-}$}

\put(24,8){ $\hat{e}_3^{-}$}

\put(-3,8){ $\hat{e}_3$}

\put(11,2){ $\widetilde{T}_{\sigma}$}
\end{picture}
}

\put(85,4){\begin{picture}(10,10)

\put(3,28){\line(1,0){15}}

\put(3,28){\line(-3,-5){5.5}}

\put(3,10){\line(1,0){15}}

\put(3,10){\line(-3,5){5.5}}

\put(18,10){\line(3,5){5.5}}

\put(18,28){\line(3,-5){5.5}}

\put(-5,24){ $\hat{e}_1$}

\put(7,30){ $\hat{e}_1^{-}$}

\put(21,23){ $\hat{e}_2$}

\put(20,12){ $\hat{e}_2^{-}$}

\put(8,5.8){ $\hat{e}_3^{-}$}

\put(-4,12){ $\hat{e}_3$}

\put(8,18){ $\curvearrowright$ }

\put(12,0){ $S_{\sigma}$}
\end{picture}
}


\end{picture}
\end{center}


The following lemma is essential in the whole paper.

\medskip

{\bf Lemma 1.1} $^{\cite {whi}}$   \ \ Every simple 3-connected
planar graph has a unique planar embedding.


\medskip
{\bf Lemma 1.2}  \ \ The minimum genus of a minor of a graph $G$ can
never be larger than $\gamma(G)$.

\medskip

{\bf Proof } \ \ Let the graph $G$ be embedded in a surface $S$,
then contracting an edge $e$ of $G$ on $S$ can obtain an embedding
of the contracted graph $G/e$ on $S$. Moreover, edge deletion can
never increase embedding genus. Thus, the lemma is obtained.
$\hspace*{\fill} \Box$

\medskip
{\bf Lemma 1.3} \ \ If an orientable surface $S$ has the form as
$(AxByCx^{-}Dy^{-}E)$, then $g(S)\geqslant1$, furthermore, the genus
of $S$ is $p(\geqslant 1)$ if, and only if, $ADCBE$ is with genus
$p-1$.

\medskip

{\bf Proof } \ \ According to the Transform 4, it is obvious.
$\hspace*{\fill} \Box$

\bigskip

 \noindent {\bf 2. The genus of the near-wheel graphs}

\bigskip

It is obvious that $W_{n}$ is 3-connected and $\gamma(W_{n})=0$. So,
according to Lemma 1.1, $W_{n}$ has an unique embedding in the
plane. We denote this unique planar embedding of $W_{n}$  by
$\mathbb{P}(W_{n})$.

\medskip

{\bf Lemma 2.1}   \ \ Let $\mathbb{P}(W_{n})$ be the planar
embedding of the wheel $W_{n}$ with $n$ spokes, $v$ be a
degree-three vertex which is not in $W_{n}$, then the genus
$\gamma(W_{n}+v)$ of the graph $W_{n}+v$ equals 0 if the three
$antennal$ $vertices$ of $W_{n}+v$ are in the same face boundary of
$\mathbb{P}(W_{n})$.

\medskip

{\bf Proof } \ \ Let $v_1$, $v_2$, $v_3$ be the three $antennal$
$vertices$ of $W_{n}+v$,  $f_1$ be the face of $\mathbb{P}(W_{n})$
with $v_1$, $v_2$, $v_3$ on it, then we can get a planar embedding
of $W_{n}+v$ by placing $v$ in the interior of $f_1$ and jointing
$vv_{i} \ (i=1,2,3)$.   $\hspace*{\fill} \Box$

\medskip

{\bf Lemma 2.2}   \ \ Let $\mathbb{P}(W_{n})$ be the planar
embedding of the wheel $W_{n}$ with $n$ spokes, $v$ be a
degree-three vertex which is not in $W_{n}$, then the genus
$\gamma(W_{n}+v)$ of the graph $W_{n}+v$ equals 1 if the following
two conditions are satisfied: (i) the three $antennal$-$vertex$ of
$W_{n}+v$ are in the boundary of two different faces of
$\mathbb{P}(W_{n})$; (ii) there is no face of $\mathbb{P}(W_{n})$
whose boundary contains all the three $antennal$-$vertex$.

\medskip

{\bf Proof } \ \ It is easy to find out that $K_{3,3}$ is a minor of
$W_{n}+v$. According to Lemma 1.2 we can get that
$\gamma(W_{n}+v)\geqslant 1$.

Let $v_1$, $v_2$, $v_3$ be the three $antennal$-$vertex$ of
$W_{n}+v$. Because the three $antennal$-$vertex$ of $W_{n}+v$ are in
the boundary of two different faces of $\mathbb{P}(W_{n})$, without
loss of generality, we may assume that $v_1$, $v_2$ are in the
boundary of $f_1$, and $v_3$ in  $f_2$, where $f_1$ and  $f_2$ are
two different faces of $\mathbb{P}(W_{n})$. Putting $v$ in the
interior of $f_1$ and joining $vv_{i} \ (i=1,2,3)$, then we will get
a torus embedding of $W_{n}+v$ if add a handle to the plane with the
edge $vv_{3}$ on it. So $\gamma(W_{n}+v)\leqslant 1$.

From the above we can get that $\gamma(W_{n}+v)\geqslant 1$ and
$\gamma(W_{n}+v)\leqslant 1$.  So $\gamma(W_{n}+v)= 1$.
$\hspace*{\fill} \Box$

\medskip

{\bf Lemma 2.3}   \ \ Let $\mathbb{P}(W_{n})$ be the planar
embedding of the wheel $W_{n}$ with $n$ spokes, $v$ be a
degree-three vertex which is not in $W_{n}$, then the genus
$\gamma(W_{n}+v)$ of the graph $W_{n}+v$ equals 1 if any pair of the
three $antennal$-$vertex$ of $W_{n}+v$ are not in a same face
boundary of $\mathbb{P}(W_{n})$.

\bigskip

{\bf Proof } \ \ It is not difficult to find out that $K_{3,3}$ is a
minor of $W_{n}+v$. According to Lemma 1.2 we can get that
$\gamma(W_{n}+v)\geqslant 1$.

\medskip

{\bf Case 1:} \ The three $antennal$-$vertex$ of $W_{n}+v$  are all
$subdividing$-$vertex$ of $W_{n}$.

\medskip

Let $v_1$, $v_2$, $v_3$ be the three $antennal$-$vertex$ of
$W_{n}+v$. For any pair of the three $antennal$-$vertex$ of
$W_{n}+v$ are not in a same face boundary of $\mathbb{P}(W_{n})$,
the vertices $v_1$, $v_2$ and $v_3$ must belong to one of the
following two subcases: (1) $v_1$, $v_2$ and $v_3$ are in three
different spokes of $W_{n}$, furthermore, any pair of these three
spokes are not in a same face boundary of $\mathbb{P}(W_{n})$; (2)
one of \{$v_1$, $v_2$, $v_3$\} is on the boundary of the unbounded
face of $\mathbb{P}(W_{n})$, and the other two are in two different
spokes of $\mathbb{P}(W_{n})$, where the two spokes are not on a
same face boundary of $\mathbb{P}(W_{n})$.

\setlength{\unitlength}{1mm}
\begin{center}
\begin{picture}(100,45)

\put(10,5){\begin{picture}(10,10)

\put(0,20){\circle*{1.5}}   

\put(0,26.5){\circle*{1.5}}  

\put(0,20){\line(0,1){14}}

\put(0,34){\circle*{1.5}}   

\put(0,34){\line(-5,-2){10}}

\put(0,34){\line(5,-2){10}}

\put(-4,29){\circle*{1.5}}

\qbezier(-4,28.3)(-6,20)(-6.8,16)

\qbezier(-4,28.3)(-2,27.5)(-0.7,26.7)

\qbezier(-4,28.3)(-5,13)(6.8,16)

\put(-13,12.3){\circle*{1.5}}  

\put(-13,12.3){\line(-1,4){1.8}}

\put(-13,12.3){\line(5,-4){6}}

\put(-6.8,16){\circle*{1.5}}

\put(0,20){\line(-5,-3){13}}   

\put(13,12.3){\circle*{1.5}}  

\put(13,12.3){\line(1,4){1.8}}

\put(13,12.3){\line(-5,-4){6}}

\put(6.8,16){\circle*{1.5}}

\put(0,20){\line(5,-3){13}}  

\put(-15,20){\circle*{1.5}}  

\put(0,20){\line(-1,0){15}}

\put(-9,21.5){\circle*{0.6}}   

\put(-8.7,23.5){\circle*{0.6}}   

\put(-8,25){\circle*{0.6}}   

\put(-10,30){\circle*{1.5}}

\put(0,20){\line(-1,1){10.2}}

\qbezier(-15,20)(-12,27.5)(-10,30)  

\begin{footnotesize} \put(-17.5,14.5){{$a_1$}}    

\put(-12,34){{$a_{m-1}$}}

\put(4,34){{$a_{m}$}}

\put(15,15){{$a_{m+p}$}}

\put(10,7){{$a_{m+p+1}$}}

\put(-13,8){{$a_{n}$}}

\put(-7.5,13){{$v_1$}}

\put(-5,30){{$v$}}

\put(1,27){{$v_2$}}

\put(7,17){{$v_3$}}   

\end{footnotesize}

\put(15,20){\circle*{1.5}}  

\put(0,20){\line(1,0){15}}

\put(9,21.5){\circle*{0.6}}   

\put(8.7,23.5){\circle*{0.6}}   

\put(8,25){\circle*{0.6}}   

\put(10,30){\circle*{1.5}}

\put(0,20){\line(1,1){10.2}}

\qbezier(15,20)(12,27.5)(10,30)  

\put(6.2,7.3){\circle*{1.5}}  

\put(0,20){\line(1,-2){6}}

\put(2,12){\circle*{0.6}}   

\put(0,11.8){\circle*{0.6}}   

\put(-2,12){\circle*{0.6}}   

\put(-6.5,7.3){\circle*{1.5}}

\put(0,20){\line(-1,-2){6.3}}

\qbezier(6.2,7.3)(0,6)(-6.5,7.3)  


\put(-12,-10){{\bf Fig.2: $W_{n}+v$}}

\end{picture}}

\put(85,6){\begin{picture}(10,10)
\put(0,20){\circle*{1.5}}   

\put(-6,22){\circle*{1.5}}

\put(-6,22){\line(-2,-1){6}}                                 

\put(0,20){\line(0,1){14}}

\put(0,34){\circle*{1.5}}   

\put(0,34){\line(-5,4){5}}

\put(0,34){\line(5,4){5}}                                 

\put(-13,12.3){\circle*{1.5}}    

\put(-13,12.3){\line(-5,-6){4}}

\put(-13,12.3){\line(-3,0){6}}

\put(0,20){\line(-5,-3){13}}

\put(-6.8,16){\circle*{1.5}}  

\put(13,12.3){\circle*{1.5}}         

\put(6.8,16){\line(1,0){6.7}}

\put(13,12.3){\line(4,-1){6.2}}

\put(13,12.3){\line(1,-3){2}}

\put(6.8,16){\circle*{1.5}}

\put(0,20){\line(5,-3){13}}         

\put(-21,16.2){\circle*{1.5}}

\put(-21,16.2){\line(1,0){13.5}}                                  

\put(-12.5,24.2){\circle*{1.5}}                                 

\put(-21,16.2){\line(-1,1){4}}

\put(-21,16.2){\line(-1,-1){4}}

\put(0,20){\line(-3,1){13}}

\put(-12.5,24.2){\line(-1,1){4}}

\put(-12.5,24.2){\line(-2,-1){5}}

\put(-1.5,30){\circle*{0.6}}   

\put(-3.3,29.7){\circle*{0.6}}   

\put(-4.8,29){\circle*{0.6}}   

\put(-8.5,31){\circle*{1.5}}

\put(0,20){\line(-3,4){8.5}}

\put(-8.5,31){\line(-6,1){7}}

\put(-8.5,31){\line(-1,4){2}}                                    

\begin{footnotesize} \put(-22.5,11){{$a_1$}}    

\put(-27.2,11.2){{$y$}}

\put(-28,19){{$x$}}

\put(-15,18.5){{$x^{-}$}}

\put(-23,22){{$a_{m}$}}

\put(-25,28){{$a_{m-1}^{-}$}}

\put(-24,33){{$a_{m-1}$}}

\put(-19,40){{$a_{m-2}^{-}$}}

\put(-6,39){{$a_{2}$}}

\put(3,40){{$a_{1}^{-}$}}

\put(11,38){{$a_{m+p}$}}

\put(16,30){{$a_{m+p-1}^{-}$}}

\put(21,22){{$a_{m+1}$}}

\put(21,16.5){{$a_{m}^{-}$}}

\put(14,15){{$y^{-}$}}

\put(20,10){{$a_{m+p}^{-}$}}

\put(15.7,5.7){{$a_{m+p+1}$}}

\put(10.3,2){{$a_{m+p+1}^{-}$}}

\put(0,-2){{$a_{m+p+2}$}}

\put(-22,6){{$a_{n}^{-}$}}

\put(-14.5,5){{$a_{n}$}}

\put(-10,0){{$a_{n-1}^{-}$}}

\put(-20,17){{$v$}}

\put(-7.5,13){{$v_1$}}

\put(-8,24){{$v_2$}}

\put(5.5,12.5){{$v_3$}}   

\end{footnotesize}

\put(15,20){\circle*{1.5}}                                    

\put(0,20){\line(1,0){15}}

\put(10.7,21.7){\circle*{0.6}}   

\put(10.3,23.7){\circle*{0.6}}   

\put(9.2,25.5){\circle*{0.6}}   

\put(10,30){\circle*{1.5}}

\put(0,20){\line(1,1){10.2}}

\put(15,20){\line(2,1){5}}

\put(15,20){\line(3,-2){5}}

\put(10,30){\line(1,6){1}}

\put(10,30){\line(1,0){5}}                                     

\put(6.2,7.3){\circle*{1.5}}  

\put(0,20){\line(1,-2){6}}

\put(2,11){\circle*{0.6}}   

\put(0,10.8){\circle*{0.6}}   

\put(-2,11){\circle*{0.6}}   

\put(-6.5,7.3){\circle*{1.5}}

\put(0,20){\line(-1,-2){6.3}}

\put(6.2,7.3){\line(-1,-2){3.3}}

\put(6.2,7.3){\line(3,-2){4}}

\put(-6.5,7.3){\line(0,-1){4}}

\put(-6.5,7.3){\line(-3,-1){4}}       


\put(-10,-11){{\bf Fig.3: $\widetilde{T}_{\sigma}$}}
\end{picture}}

\end{picture}
\end{center}

\bigskip


In the first subcase, the graph $W_{n}+v$ and one of its joint-tree
are shown in Fig.2 and Fig.3 respectively, where we denoted the edge
($v$, $v_2$) by $x$, and ($v$, $v_3$) by $y$. In Fig.2, the edges of
the $n$-cycle in $W_{n}$, according to the clockwise rotation, are
denoted by $a_1$, $a_2$, $\dots$, $a_{n}$. The surface associated
with the joint-tree in Fig.3 is
\begin{eqnarray}
S & = & a_1yxx^{-}a_{m}a_{m-1}^{-}a_{m-1}a_{m-2}^{-}a_{m-2} \dots
a_{2}^{-}a_{2}a_{1}^{-}a_{m+p}a_{m+p-1}^{-}a_{m+p-1}a_{m+p-2}^{-}
\nonumber\\
& & a_{m+p-2} \dots
a_{m+1}^{-}a_{m+1}a_{m}^{-}y^{-}a_{m+p}^{-}a_{m+p+1}a_{m+p+1}^{-}a_{m+p+2}a_{m+p+2}^{-}
\dots
a_{n}a_{n}^{-}  \nonumber \\
& \sim &
a_1ya_{m}a_{1}^{-}a_{m+p}a_{m}^{-}y^{-}a_{m+p}^{-} \nonumber \\
& \sim &
a_{m+p}a_{m}^{-}a_{m}a_{m+p}^{-}a_1ya_{1}^{-}y^{-} \nonumber \\
& \sim & a_1ya_{1}^{-}y^{-} \nonumber
\end{eqnarray}

Obviously, $g(S)=1$. So $\gamma(W_{n}+v)\leqslant 1$. On the other
hand $\gamma(W_{n}+v)\geqslant 1$. Therefore, in the first subcase,
$\gamma(W_{n}+v)= 1$.

\medskip

In the second subcase, the graph $W_{n}+v$ and one of its joint-tree
are shown in Fig.4 and Fig.5 respectively, where we denoted the edge
($v$, $v_2$) by $x$, and ($v$, $v_3$) by $y$. In Fig.4, the edges of
the $n$-cycle in $W_{n}$, according to the clockwise rotation, are
denoted by $a_1$, $a_2$, $\dots$, $a_{m-1}$, $b$, $a_{m}$, $\dots$,
$a_{n}$. The surface associated with the joint-tree in Fig.5 is
\begin{eqnarray}
S & = & a_1yxx^{-}a_{m}a_{m-1}^{-}a_{m-1}a_{m-2}^{-}a_{m-2} \dots
a_{2}^{-}a_{2}a_{1}^{-}a_{m+p}a_{m+p-1}^{-}a_{m+p-1}a_{m+p-2}^{-}
\nonumber\\
& & a_{m+p-2} \dots
a_{m+1}^{-}a_{m+1}a_{m}^{-}y^{-}a_{m+p}^{-}a_{m+p+1}a_{m+p+1}^{-}a_{m+p+2}a_{m+p+2}^{-}
\dots
a_{n}a_{n}^{-}  \nonumber \\
& \sim &
a_1ya_{m}a_{1}^{-}a_{m+p}a_{m}^{-}y^{-}a_{m+p}^{-} \nonumber \\
& \sim &
a_{m+p}a_{m}^{-}a_{m}a_{m+p}^{-}a_1ya_{1}^{-}y^{-} \nonumber \\
& \sim & a_1ya_{1}^{-}y^{-} \nonumber
\end{eqnarray}

Obviously, $g(S)=1$. So $\gamma(W_{n}+v)\leqslant 1$. On the other
hand $\gamma(W_{n}+v)\geqslant 1$. Therefore, in the second subcase,
$\gamma(W_{n}+v)= 1$.

\medskip


\setlength{\unitlength}{1mm}
\begin{center}
\begin{picture}(100,45)

\put(10,5){\begin{picture}(10,10)

\put(0,20){\circle*{1.5}}   

\put(0,20){\line(0,1){14}}

\put(0,34){\circle*{1.5}}   

\put(0,34){\line(5,-2){10}}

\put(-18,29){\circle*{1.5}}

\qbezier(-18,29)(-13,18)(-6.8,16)

\qbezier(-18,29)(-15,36)(-8,32)

\qbezier(-18,29)(-2,29)(6.8,16)      

\put(-13,12.3){\circle*{1.5}}  

\put(-13,12.3){\line(-1,4){1.8}}

\put(-13,12.3){\line(5,-4){6}}

\put(-6.8,16){\circle*{1.5}}  

\put(0,20){\line(-5,-3){13}}   

\put(13,12.3){\circle*{1.5}}  

\put(13,12.3){\line(1,4){1.8}}

\put(13,12.3){\line(-5,-4){6}}

\put(6.8,16){\circle*{1.5}}

\put(0,20){\line(5,-3){13}}  

\put(-15,20){\circle*{1.5}}  

\put(0,20){\line(-1,0){15}}

\put(-11.2,21){\circle*{0.6}}   

\put(-11.1,22.5){\circle*{0.6}}   

\put(-10.6,23.9){\circle*{0.6}}   

\put(-8,32){\circle*{1.5}}  

\put(-13,27){\circle*{1.5}}

\put(0,20){\line(-2,1){13}}

\qbezier(-7.5,32.3)(-4,34)(0,34)

\qbezier(-12.8,27.5)(-11,30)(-8.4,31.5)

\qbezier(-15,20)(-12.8,27.5)(-13,27)  

\begin{footnotesize} \put(-17.5,14.5){{$a_1$}}    

\put(-7,35){{$a_{m}$}}

\put(3,34){{$a_{m+1}$}}

\put(15,15){{$a_{m+p}$}}

\put(10,7){{$a_{m+p+1}$}}

\put(-13,8){{$a_{n}$}}

\put(-12.8,30){{$b$}}

\put(-7.5,13){{$v_1$}}

\put(-22,29){{$v$}}

\put(-7,29){{$v_2$}}

\put(7,17){{$v_3$}}   

\end{footnotesize}

\put(15,20){\circle*{1.5}}  

\put(0,20){\line(1,0){15}}

\put(9,21.5){\circle*{0.6}}   

\put(8.7,23.5){\circle*{0.6}}   

\put(8,25){\circle*{0.6}}   

\put(10,30){\circle*{1.5}}

\put(0,20){\line(1,1){10.2}}

\qbezier(15,20)(12,27.5)(10,30)  

\put(6.2,7.3){\circle*{1.5}}  

\put(0,20){\line(1,-2){6}}

\put(2,12){\circle*{0.6}}   

\put(0,11.8){\circle*{0.6}}   

\put(-2,12){\circle*{0.6}}   

\put(-6.5,7.3){\circle*{1.5}}

\put(0,20){\line(-1,-2){6.3}}

\qbezier(6.2,7.3)(0,6)(-6.5,7.3)  


\put(-12,-10){{\bf Fig.4: $W_{n}+v$}}

\end{picture}}


\put(85,6){\begin{picture}(10,10)
\put(0,20){\circle*{1.5}}   

\put(0,20){\line(0,1){14}}

\put(0,34){\circle*{1.5}}   

\put(0,34){\line(-5,4){5}}

\put(0,34){\line(5,4){5}}                                 

\put(-13,12.3){\circle*{1.5}}    

\put(-13,12.3){\line(-5,-6){4}}

\put(-13,12.3){\line(-3,0){6}}

\put(0,20){\line(-5,-3){13}}

\put(-6.8,16){\circle*{1.5}}  

\put(13,12.3){\circle*{1.5}}         

\put(6.8,16){\line(1,0){6.7}}

\put(13,12.3){\line(4,-1){6.2}}

\put(13,12.3){\line(1,-3){2}}

\put(6.8,16){\circle*{1.5}}

\put(0,20){\line(5,-3){13}}         

\put(-21,16.2){\circle*{1.5}}

\put(-21,16.2){\line(1,0){13.5}}                                 

\put(-12.5,24.2){\circle*{1.5}}                                 

\put(-21,16.2){\line(-5,2){6}}

\put(-21,16.2){\line(-6,-5){5}}                 

\put(0,20){\line(-3,1){13}}

\put(-13,24.2){\line(-1,0){7.5}}

\put(-12.5,24.2){\line(-1,1){4}}      

\put(-1.5,30){\circle*{0.6}}   

\put(-3.3,29.7){\circle*{0.6}}   

\put(-4.8,29){\circle*{0.6}}   

\put(-20.5,24.2){\line(-5,2){7}}

\put(-20.5,24.2){\circle*{1.5}}                                   

\put(-8.5,31){\circle*{1.5}}

\put(-21,24.2){\line(-3,-1){7}}    

\put(0,20){\line(-3,4){8.5}}

\put(-8.5,31){\line(-3,1){7}}

\put(-8.5,31){\line(-1,4){2}}                                     

\begin{footnotesize} \put(-23,11){{$a_1$}}    

\put(-30,12){{$y$}}

\put(-31,16.5){{$x$}}

\put(-32.5,21){{$x^{-}$}}

\put(-22,21.5){{$v_{2}$}}

\put(-17,21){{$b$}}

\put(-33.5,26.5){{$a_{m}$}}

\put(-25,29){{$a_{m-1}^{-}$}}

\put(-24,35){{$a_{m-1}$}}

\put(-19,40){{$a_{m-2}^{-}$}}

\put(-6,39){{$a_{2}$}}

\put(3,40){{$a_{1}^{-}$}}

\put(11,38){{$a_{m+p}$}}

\put(16,30){{$a_{m+p-1}^{-}$}}

\put(21,22){{$a_{m+1}$}}

\put(21,16.5){{$a_{m}^{-}$}}

\put(14,15){{$y^{-}$}}

\put(20,10){{$a_{m+p}^{-}$}}

\put(15.7,5.7){{$a_{m+p+1}$}}

\put(10.3,2){{$a_{m+p+1}^{-}$}}

\put(0,-2){{$a_{m+p+2}$}}

\put(-22,6){{$a_{n}^{-}$}}

\put(-14.5,5){{$a_{n}$}}

\put(-10,0){{$a_{n-1}^{-}$}}

\put(-20,17){{$v$}}

\put(-7.5,13){{$v_1$}}

\put(5.5,12.5){{$v_3$}}   

\end{footnotesize}

\put(15,20){\circle*{1.5}}                                    

\put(0,20){\line(1,0){15}}

\put(10.7,21.7){\circle*{0.6}}   

\put(10.3,23.7){\circle*{0.6}}   

\put(9.2,25.5){\circle*{0.6}}   

\put(10,30){\circle*{1.5}}

\put(0,20){\line(1,1){10.2}}

\put(15,20){\line(2,1){5}}

\put(15,20){\line(3,-2){5}}

\put(10,30){\line(1,6){1}}

\put(10,30){\line(1,0){5}}                                     

\put(6.2,7.3){\circle*{1.5}}  

\put(0,20){\line(1,-2){6}}

\put(2,11){\circle*{0.6}}   

\put(0,10.8){\circle*{0.6}}   

\put(-2,11){\circle*{0.6}}   

\put(-6.5,7.3){\circle*{1.5}}

\put(0,20){\line(-1,-2){6.3}}

\put(6.2,7.3){\line(-1,-2){3.3}}

\put(6.2,7.3){\line(3,-2){4}}

\put(-6.5,7.3){\line(0,-1){4}}

\put(-6.5,7.3){\line(-3,-1){4}}      


\put(-10,-11){{\bf Fig.5: $\widetilde{T}_{\sigma}$}}
\end{picture}}

\end{picture}
\end{center}
\bigskip

According to the above, we can get that, in the Case 1,
$\gamma(W_{n}+v)= 1$.

\medskip

{\bf Case 2:} \ The three $antennal$-$vertex$ of $W_{n}+v$ consist
of both $subdividing$-$vertex$ of $W_{n}$ and vertices which belong
to $V(W_{n})$.

\medskip

Because any pair of the three $antennal$-$vertex$ of $W_{n}+v$ are
not in a same face boundary of $\mathbb{P}(W_{n})$, among these
three $antennal$ $vertices$, there is one and only one vertex
belongs to $V(W_{n})$, and the other two are both
$subdividing$-$vertex$ of $W_{n}$. It is not difficult to find out
that the graph $W_{n}+v$ in Case 2 is minor of the graph $W_{n}+v$
in Case 1. So, according to Lemma 1.2 we can get that, in Case 2,
$\gamma(W_{n}+v)\leqslant1$. On the other hand, we can get that
$\gamma(W_{n}+v)\geqslant1$  because $K_{3,3}$ is a minor of
$W_{n}+v$. So, in the Case 2, $\gamma(W_{n}+v)=1$.

According to the Case 1 and Case 2 we can get the Lemma 2.3.
$\hspace*{\fill} \Box$

\medskip

The following theorem can be easily  obtained from Lemma 2.1, Lemma
2.2 and Lemma 2.3.

\medskip

{\bf Theorem A}   \  Let $\mathbb{P}(W_{n})$ be the planar embedding
of the wheel $W_{n}$ with $n$ spokes, $v$ be a degree-three vertex
which is not in $W_{n}$, then the genus $\gamma(W_{n}+v)$ of the
graph $W_{n}+v$ equals 0 if the three $antennal$-$vertex$ of
$W_{n}+v$ are in the same face boundary of $\mathbb{P}(W_{n})$,
otherwise, $\gamma(W_{n}+v)= 1$.

\medskip

{\bf Remark}   \ \ (i) From theorem A we can get that there are many
planar or toroidal graphs whose genus range can be arbitrarily
large; (ii) How does the genus of a cubic planar graph $G$ change if
we add a degree-three vertex $v$, which is not in $V(G)$, to $G$? We
believe its genus to be 0 or 1. So, the proof or disproof of the
result will be interesting.

\bigskip

 \noindent {\bf 3. Lower bound on the maximum genus of graphs}

\bigskip

A set $J\subseteq V(G)$ is called a $non$-$separating$ $independent$
$set$ of a connected graph $G$ if $J$ is an independent set of $G$
and $G-J$ is connected. In 1997, through the independent set of a
graph, Huang and Liu$^{\cite{hua}}$ studied the maximum genus of
cubic graphs, and obtained the following result.

\medskip

{\bf Lemma 3.1 }$^{\cite{hua}}$  \ The maximum genus of a cubic
graph $G$ equals the cardinality of the maximum non-separating
independent set of $G$.
\medskip

\hspace{-8mm}  But for general graphs that is not necessary cubic,
there is no result concerning the maximum genus which is
characterized by the independent set of the graph. In the following,
we will provide a lower bound of the maximum genus, which is
characterized via the independent set, for general graphs.
Furthermore, there are examples shown that the bound may be tight,
and, in some sense,  may be better than the result obtained by Li
and Liu$^{\cite{li}}$, and the result obtained by Z. Ouyang
$etc.^{\cite{ou}}$.

\bigskip

{\bf Theorem B}   \  Let $G$ be a connected graph whose minimum
degree is at leas 3. If $A=\{v_1, v_2,  \dots, v_{m}\}$  is an
independent set such that $G-A$ is connected,then
\begin{displaymath}
 \gamma_{M}( G ) \geqslant \frac{1}{2}\sum_{i=1}^{m}\big(d(v_{i})-\varepsilon_{i}\big)+  \gamma_{M}\big(G-\{v_1, v_2,  \dots,
 v_{m}\}\big),
\end{displaymath}
where for each index $i (1\leqslant i \leqslant m)$,
$\varepsilon_{i}=1$ if $d(v_{i})\equiv 1 (mod \ 2)$ and
$\varepsilon_{i}=2$ otherwise.

\medskip

{\bf Proof}   \ Without loss of generality, let $H$ be the graph
obtained from $G$ by successively deleting $v_1, v_2,  \dots, v_{m}$
from $G$, and $\psi(H)$ be a maximum genus embedding of $H$. We
first add the vertex $v_{m}$ to $H$.

\medskip

\textbf{Case 1:} \ $d_{G}(v_{m})\equiv 1 \ (mod \ 2)$.

\medskip

Without loss of generality, let $d_{G}(v_{m})=2i+1$, and  $x_1, x_2,
\dots, x_{2i+1}$ be the $2i+1$ neighbors of $v_{m}$ in $G$.
According to the $2i+1$ neighbors of $v_{m}$ are in the same face
boundary of $\psi(H)$ or not, we will discuss in the following two
subcases.

\medskip

\textsf{Subcase 1.1:} \ All the neighbors of $v_{m}$ are in the same
face boundary of $\psi(H)$.

\medskip

Let $f_0$, which is bounded by $B_0$, be the face of $\psi(H)$ that
$x_1, x_2, \dots, x_{2i+1}$ are on the boundary of it. Firstly, we
put $v_{m}$ in $f_0$ and connect each of \{$x_1, x_2, x_3$\} to
$v_{m}$, and denote this resulting graph by $H_1$. Through the
manner depicted by Fig.7, where each $vertex$-$rotation$ is the same
with that of $\psi(H)$ except $v_{m}$, we can get an embedding
$\psi(H_1)$ of $H_1$ such that its face number is the same with that
of $\psi(H)$. From the equation $V-E+F=2-2g$, it can be easily
deduced that the maximum genus of $H_1$ is at least one more than
that of $H$.

Now connect each of  \{$x_4, x_5$\} to $v_{m}$, and denote the
resulting graph by $H_2$. Through the manner depicted by Fig.8, we
can get an embedding $\psi(H_2)$ of $H_2$, which has the same face
number with that of $\psi(H)$. From the equation $V-E+F=2-2g$, it
can be easily deduced that the maximum genus of $H_2$ is at least
two more than that of $H$.

\medskip


\setlength{\unitlength}{1mm}
\begin{center}
\begin{picture}(100,41)

\put(-5,5) {\begin{picture}(10,10)

\put(0,34){\circle*{1.5}}   

\put(-17.5,21){\circle*{0.6}}   

\put(-18,17){\circle*{0.6}}   

\put(-16.5,13){\circle*{0.6}}   

\put(-13,27){\circle*{1.5}}   

\qbezier(-12.8,27.5)(-8,34)(0,34)  

\begin{footnotesize}

\put(0,36){{$x_1$}}

\put(17,20){{$x_2$}}

\put(15,10.5){{$x_3$}}

\put(-1,2.5){{$x_4$}}

\put(-14,6){{$x_5$}}

\put(-22,28){{$x_{2i+1}$}}

\put(-1,16) {{$f_0$}}

\end{footnotesize}

\put(15,20){\circle*{1.5}}  

\put(13,12.3){\circle*{1.5}}  

\put(-11,9){\circle*{1.5}}  

\put(0,6){\circle*{1.5}}  

\qbezier(15,20)(13,33)(0,34)  

\qbezier(15,20)(15,16)(13,12.3)   

\qbezier(13,12.3)(9,6)(0,6)  

\qbezier(-11,9)(-8,6)(0,6)  

\qbezier(-10,21)(-3,34)(7,23)   

\qbezier(-11,9)(-18,17)(-13,27)  

\thicklines

\put(4.1,26){\vector(1,-1){3.5}}

\thinlines

\put(-9,-8){{\bf Fig.6: $B_0$}}

\end{picture} }


\put(50,5){\begin{picture}(10,10)

\put(0,20){\circle*{1.5}}   

\put(0,20){\line(0,1){14}}  

\put(0,34){\circle*{1.5}}   

\put(-17.5,21){\circle*{0.6}}   

\put(-18,17){\circle*{0.6}}   

\put(-16.5,13){\circle*{0.6}}   

\put(-13,27){\circle*{1.5}}   

\qbezier(-12.8,27.5)(-8,34)(0,34)  

\begin{footnotesize}

\put(0,36){{$x_1$}}

\put(17,20){{$x_2$}}

\put(15,10.5){{$x_3$}}

\put(-1,2.5){{$x_4$}}

\put(-14,6){{$x_5$}}

\put(-22,28){{$x_{2i+1}$}}

\put(-6,21){{$v_{m}$}}

\end{footnotesize}

\put(15,20){\circle*{1.5}}  

\put(13,12.3){\circle*{1.5}}  

\put(-11,9){\circle*{1.5}}  

\put(0,6){\circle*{1.5}}  

\put(0,20){\line(1,0){15}}  

\qbezier(15,20)(13,33)(0,34)  

\qbezier(15,20)(15,16)(13,12.3)   

\qbezier(13,12.3)(9,6)(0,6)  

\qbezier(-11,9)(-8,6)(0,6)  

\qbezier(-11,9)(-18,17)(-13,27)  

\qbezier(0,20)(8,34)(13,12.3)   


\put(-11,-8){{\bf Fig.7: $\psi(H_1)$}}

\end{picture}}


\put(105,5){\begin{picture}(10,10)

\put(0,20){\circle*{1.5}}   

\put(0,20){\line(0,1){14}}  

\put(0,34){\circle*{1.5}}   

\put(-17.5,21){\circle*{0.6}}   

\put(-18,17){\circle*{0.6}}   

\put(-16.5,13){\circle*{0.6}}   

\put(-13,27){\circle*{1.5}}   

\qbezier(-12.8,27.5)(-8,34)(0,34)  

\begin{footnotesize}

\put(0,36){{$x_1$}}

\put(17,20){{$x_2$}}

\put(15,10.5){{$x_3$}}

\put(-1,2.5){{$x_4$}}

\put(-14,6){{$x_5$}}

\put(-22,28){{$x_{2i+1}$}}

\put(-6,21){{$v_{m}$}}

\end{footnotesize}

\put(15,20){\circle*{1.5}}  

\put(13,12.3){\circle*{1.5}}  

\put(-11,9){\circle*{1.5}}  

\put(0,6){\circle*{1.5}}  

\put(0,20){\line(0,-1){14}}  

\put(0,20){\line(1,0){15}}  

\qbezier(15,20)(13,33)(0,34)  

\qbezier(15,20)(15,16)(13,12.3)   

\qbezier(13,12.3)(9,6)(0,6)  

\qbezier(-11,9)(-8,6)(0,6)  

\qbezier(0,20)(8,34)(13,12.3)   

\qbezier(-11,9)(-18,17)(-13,27)  

\qbezier(-11,9)(19,11)(0,20)  


\put(-11,-8){{\bf Fig.8: $\psi(H_2)$}}

\end{picture}}

\end{picture}
\end{center}

\medskip

Similar to the manner of connecting \{$x_4, x_5$\} to $v_{m}$, we
can connect \{$x_6, x_7$\}, \dots,  \{$x_{2i}, x_{2i+1}$\} to
$v_{m}$. Eventually, we will get an embedding of $H+v_{m}$. It can
be easily deduced that the maximum genus of $H+v_{m}$ is at least
$\frac{1}{2}\big(d(v_{m})-1\big)+ \gamma_{M}\big(G-\{v_1, v_2,
\dots, v_{m}\}\big)$.

\medskip

\textsf{Subcase 1.2:} \ There is no face boundary of $\psi(H)$
containing all the neighbors of $v_{m}$.

\medskip

First, add $v_{m}$ to $H$ and connect each of $\{x_1, x_2, x_3\}$ to
$v_{m}$. The resulting graph is denoted by $H_1$. If $x_1, x_2, x_3$
are in two different face boundaries of $\psi(H)$, say $f_1$ and
$f_2$, then via the manner depicted by the left part of Fig.9, we
can get an embedding $\psi(H_1)$ of $H_1$ whose face number is the
same with that of $\psi(H)$. If $x_1, x_2, x_3$ are in three
different face boundaries of $\psi(H)$, say $f_1$, $f_2$, and $f_3$,
then through the manner depicted by the right part of Fig.9, we can
get an embedding $\psi(H_1)$ of $H_1$ whose face number is two less
than that of $\psi(H)$. From the equation $V-E+F=2-2g$, it can be
easily deduced that the maximum genus of $H_1$ is at least one more
than that of $H$.

\medskip


\setlength{\unitlength}{1mm}
\begin{center}
\begin{picture}(100,32)

\put(-6,5){\begin{picture}(10,10)

\put(0,14){\circle{14}}   

\put(25,14){\circle{14}}

\qbezier(0,13)(-17,23)(32,13)   

\put(2,20.5){\circle*{1.5}}   

\put(0,13){\line(1,4){2}}   

\put(0,13){\line(1,-3){2}}   

\put(2,7.3){\circle*{1.5}}   

\put(32,13){\circle*{1.5}}   

\put(0,13){\circle*{1.5}}   

\begin{footnotesize}

\put(3,5){{$x_1$}}

\put(2,22){{$x_2$}}

\put(24,0){{$f_{2}$}}

\put(34,12){{$x_3$}}

\put(1.3,12){{$v_{m}$}}

\put(-1,0){{$f_{1}$}}

\end{footnotesize}


\put(48,-8){{\bf Fig.9}}

\end{picture}}


\put(78,5){\begin{picture}(10,10)

\put(0,6){\circle{14}}   

\put(25,6){\circle{14}}

\put(12,22){\circle{14}}

\put(13,22){\circle*{1.5}}   

\qbezier(13,22)(11,11)(-5.6,2)  

\put(13,22){\line(1,0){6}}   

\put(19,22){\circle*{1.5}}   

\qbezier(13,22)(-4,30)(27,-0.6)  

\put(-5.6,2){\circle*{1.5}}   

\put(27,-0.6){\circle*{1.5}}   

\begin{footnotesize}

\put(12,24){{$v_{m}$}}

\put(20.7,22){{$x_{1}$}}

\put(28,-3){{$x_{2}$}}

\put(-8,-1){{$x_{3}$}}

\put(-2,22){{$f_{1}$}}

\put(-15,4){{$f_{2}$}}

\put(33,4){{$f_{3}$}}

\end{footnotesize}

\thinlines

\end{picture}}
\end{picture}
\end{center}

\medskip

Similarly, connect $\{x_4, x_5\}$, \dots, $\{x_{2i}, x_{2i+1}\}$ to
$v_{m}$. Eventually, we will get an embedding of $H+v_{m}$, and it
can be easily deduced that the maximum genus of $H+v_{m}$ is at
least $\frac{1}{2} \big( d(v_{m})-1 \big)+ \gamma_{M}\big(G-\{v_1,
v_2, \dots,
 v_{m}\}\big)$.

From Subcase 1.1 and Subcase 1.2 we can get that if $d_{G}(v_{m})=1
\ (mod \ 2)$, then $\gamma_{M}(H+v_{m}) \geqslant
\frac{1}{2}(d(v_{m})-1)+ \gamma_{M}(G-\{v_1, v_2,  \dots,
 v_{m}\})$.

\medskip

\textbf{Case 2:} \ $d_{G}(v_{m})\equiv0 \ (mod \ 2)$.

\medskip

Similar to that of Case 1, we can get that if $d_{G}(v_{m})\equiv0 \
(mod \ 2)$, then $\gamma_{M}(H+v_{m}) \geqslant
\frac{1}{2}(d(v_{m})-2)+ \gamma_{M}(G-\{v_1, v_2,  \dots,
 v_{m}\})$.

From Case 1 and Case 2 we can get that $\gamma_{M}(H+v_{m})
\geqslant \frac{1}{2}(d(v_{m})-\varepsilon_{i})+ \gamma_{M}(H)$,
where $\varepsilon_{i}=1$ if $d(v_{i})\equiv 1 (mod \ 2)$ and
$\varepsilon_{i}=2$ otherwise.

\medskip

Similarly to that of $v_{m}$, we can add $v_{m-1}$, $v_{m-2}$,
\dots, $v_{1}$, one by one, to $H+v_{m}$. Eventually  we will get an
embedding of $G$, and it is not hard to obtain that the maximum
genus of $G$ is at least
$\frac{1}{2}\sum_{i=1}^{m}(d(v_{i})-\varepsilon_{i})+
\gamma_{M}(G-\{v_1, v_2,  \dots,
 v_{m}\})$, where for each index $i (1\leqslant i \leqslant m)$,
$\varepsilon_{i}=1$ if $d(v_{i})\equiv 1 (mod \ 2)$ and
$\varepsilon_{i}=2$ otherwise. $\hspace*{\fill} \Box$

\bigskip

Noticing that the upper embeddability of a graph would not be
changed if adding an odd vertex to it, we can get the following
theorem whose proof is similar to that of Theorem B.

\bigskip

{\bf Theorem C}   \ Let $G$ be a connected graph and $A_1, A_2,
\dots A_{s}$ be a sequence of disjoint independent vertex sets which
satisfy: (i) $G_0=G$, $G_{i}=G_{i-1}-A_{i}$ is connected $(i=1, 2,
\dots, s)$; (ii) each vertex of $A_{i}$ $(i=1, 2, \dots, s)$ is an
odd vertex in $G_{i-1}$. Then for $i=0, 1, \dots, s-1$,
\begin{displaymath}
 \gamma_{M}( G_{i} ) \geqslant \frac{1}{2}\sum_{v \in A_{i+1}} \big(d_{G_{i}}(v)-1 \big )+
 \gamma_{M}(G_{i+1}).
\end{displaymath}
In particular, if one of the graph sequence $G_1, G_2, \dots, G_{s}$
is upper embeddable, then $G$ is upper embeddable.

\bigskip

{\bf Remark}   \  \ In 2000, through the girth $g$ and connectivity
of graphs, D. Li and Y. Liu$^{\cite{li}}$ obtained the lower bound
of the maximum genus of graphs, which is displayed by the following
table, where the first row and the first column represents the girth
and connectivity respectively.

\begin{center}
\scriptsize
\begin{tabular*}{154mm}[9mm]{|l|lllllllll|}

\hline

 & $g$=3  &  $g$=4  &  $g$=5  &  $g$=6  & $g$= 7 & $g$= 8 &  $g$=9 &  $g$=10 &  $g$=12  \\

\hline

1 &  $\frac{ \beta(G) + 2 }{4}$  &  $\frac{ \beta(G) + 2 }{3}$  &
$\frac{ 2\beta(G) + 2 }{5}$ & $\frac{ 3\beta(G) + 2 }{7}$ & $\frac{
5\beta(G) + 2 }{11}$ & $\frac{ 7\beta(G) + 2 }{15}$
& $\frac{ 14\beta(G) + 2 }{29}$ & $\frac{ 17\beta(G) + 2 }{35}$ & $\frac{ 31\beta(G) + 2 }{63}$ \\

2  &  $\frac{ \beta(G) + 2 }{3}$  &  $\frac{ \beta(G) + 2 }{3}$  &
$\frac{ 2\beta(G) + 3 }{5}$ & $\frac{ 3\beta(G) + 4 }{7}$ & $\frac{
6\beta(G) + 7 }{13}$ & $\frac{ 7\beta(G) + 8 }{15}$
& $\frac{ 14\beta(G) + 15 }{29}$ & $\frac{ 17\beta(G) + 18 }{35}$ & $\frac{ 31\beta(G) + 32 }{63}$ \\

3 &  $\frac{ \beta(G) + 2 }{3}$  &  $\frac{ 3\beta(G) + 4 }{7}$  &
$\frac{ 5\beta(G) + 6 }{11}$ & $\frac{ 7\beta(G) + 8 }{15}$ &
$\frac{ 11\beta(G) + 12 }{23}$ & $\frac{ 15\beta(G) + 16 }{31}$
& $\frac{ 29\beta(G) + 30 }{59}$ & $\frac{ 35\beta(G) + 36 }{71}$ & $\frac{ 63\beta(G) + 64 }{127}$ \\

\hline
\end{tabular*}
\end{center}

\medskip
Ten years later,  Z. Ouyang, J. Wang and Y. Huang$^{\cite{ou}}$
studied this parameter too, and obtained that: Let $G$ be a
$k$-$edge$-$connected$ (or $k$-$connected$) simple graph with
minimum degree $\delta$ and girth $g$. Then
 $\gamma_{M}(G)\geqslant min\{f_{k}(\delta,g)(\beta(G)+1),
\lfloor\frac{\beta(G)}{2}\rfloor \}$ for $k=1,2,3,$ where

\begin{center}
\begin{tabular*}{153.5mm}[20mm]{|l|l|l|l|}

\hline

$\delta$ \ \ & $f_{1}(\delta,g)$   \ \ \ & $f_{2}(\delta,g)$   \ \ \ & $f_{3}(\delta,g)$   \\

\hline

$\delta=3$  &  $\frac{1}{4}$  & $\frac{1}{3}$  &
$\frac{1}{2}(1-\frac{1}{4\lceil\frac{(\delta-2)(\delta+g-3)-3}{4}\rceil+3}$) \\

$\delta\geqslant4$  &
$\frac{1}{2}(1-\frac{3}{4\lceil\frac{(\delta-2)(\delta+g-3)-3}{4}\rceil+1}$)
&
$\frac{1}{2}(1-\frac{1}{2\lceil\frac{(\delta-2)(\delta+g-3)-3}{4}\rceil+1}$)
&
$\frac{1}{2}(1-\frac{1}{4\lceil\frac{(\delta-2)(\delta+g-3)-3}{4}\rceil+3}$) \\

\hline
\end{tabular*}.
\end{center}

\medskip

There are many examples showing that the lower bound in Theorem B
may be best possible. Furthermore, it may be better than the result
obtained by Li and Liu$^{\cite{li}}$ and the result of Z. Ouyang
$etc.^{\cite{ou}}$. The following are two examples with girth 3 and
connectivity 2, and girth 4 and connectivity 3 respectively.

\medskip

\setlength{\unitlength}{1mm}
\begin{center}
\begin{picture}(100,27)

\put(-14,6){\begin{picture}(10,10)

\put(0,20){\circle*{1.8}}   

\put(12,20){\circle*{1.5}}   

\put(0,5){\circle*{1.5}}   

\put(12,5){\circle*{1.5}}   

\put(0,20){\line(1,0){12}}   

\put(0,20){\line(0,-1){15}}   

\put(0,20){\line(4,-5){12}}   

\put(12,5){\line(-1,0){12}}   

\put(12,5){\line(1,0){20}}   

\put(12,5){\line(0,1){15}}   

\put(0,5){\line(4,5){12}}   

\put(22,12){\circle*{1.5}}   

\put(32,20){\circle*{1.5}}   

\put(44,20){\circle*{1.8}}   

\put(32,5){\circle*{1.5}}   

\put(44,5){\circle*{1.5}}   

\put(32,20){\line(1,0){12}}   

\put(32,20){\line(0,-1){15}}   

\put(32,20){\line(4,-5){12}}   

\put(44,5){\line(-1,0){12}}   

\put(44,5){\line(0,1){15}}   

\put(32,5){\line(4,5){12}}   

\put(22,12){\line(5,4){10}}   

\put(22,12){\line(-5,4){10}}   

\put(22,12){\line(4,-3){10}}   

\begin{footnotesize}

\put(44,22){{$v_2$}}

\put(-1,22){{$v_1$}}

\put(3,-5){{\textsf{\bf} girth 3 and 2-connected}}

\end{footnotesize}


\put(57,-8){{\bf Fig.10}}

\end{picture}}


\put(80,5){\begin{picture}(10,10)

\put(0,18){\circle*{1.8}}   

\put(0,25){\circle*{1.5}}   

\put(-10,18){\circle*{1.5}}   

\put(-10,25){\circle*{1.8}}   

\put(-10,25){\line(1,0){10}}   

\put(-10,25){\line(0,-1){7}}   

\qbezier(-10,25)(-20,14)(-10,3)

\put(-10,10){\line(1,0){10}}   

\put(-10,10){\line(0,1){10}}   

\put(-10,10){\line(0,-1){7}}   

\put(0,18){\line(-1,0){10}}   

\put(0,18){\line(0,1){7}}   

\put(0,3){\line(0,1){7}}   

\put(0,3){\line(-1,0){10}}   

\put(0,3){\line(1,0){20}}   

\put(0,10){\circle*{1.8}}  

\put(0,3){\circle*{1.5}}   

\put(-10,10){\circle*{1.5}}   

\put(-10,3){\circle*{1.5}}   

\put(10,14){\circle*{1.5}}   

\put(10,14){\line(5,2){10}}   

\put(10,14){\line(5,-2){10}}   

\put(10,14){\line(-5,2){10}}   

\put(10,14){\line(-5,-2){10}}   

\put(20,18){\circle*{1.8}}   

\put(20,25){\circle*{1.5}}   

\put(30,18){\circle*{1.5}}   

\put(30,25){\circle*{1.5}}   

\put(20,10){\circle*{1.5}}   

\put(20,3){\circle*{1.5}}   

\put(30,10){\circle*{1.8}}   

\put(30,3){\circle*{1.5}}   

\put(20,25){\line(1,0){10}}   

\put(20,25){\line(-1,0){20}}   

\put(20,25){\line(0,-1){7}}   

\put(30,18){\line(0,1){7}}   

\put(30,18){\line(0,-1){8}}   

\put(30,18){\line(-1,0){10}}   

\put(30,3){\line(-1,0){10}}   

\put(30,3){\line(0,1){7}}   

\put(20,10){\line(1,0){10}}   

\put(20,10){\line(0,-1){7}}   

\qbezier(30,25)(40,14)(30,3)

\begin{footnotesize}

\put(1,19){{$v_{2}$}}

\put(-11,27){{$v_{1}$}}

\put(-2,12){{$v_{3}$}}

\put(15.5,19.5){{$v_{4}$}}

\put(26,12){{$v_{5}$}}

\put(-10,-4){{\textsf{\bf} girth 4 and 3-connected}}

\end{footnotesize}

\thinlines

\end{picture}}
\end{picture}
\end{center}

\medskip

In the graph $G$ depicted in the left of Fig.10, let $A=\{v_1,
v_2\}$. Then

\vspace{-5mm}
\begin{eqnarray*}
\frac{1}{2}\sum_{i=1}^{m}\big(d(v_{i})-\varepsilon_{i}\big)+
\gamma_{M}\big(G-\{v_1, v_2,  \dots,
 v_{m}\}\big) \\
\lefteqn{ = \frac{1}{2}\big( (3-1)+(3-1)  \big) + \gamma_{M} \big(
G-\{v_1, v_2 \}  \big)   }  \hspace*{86mm} \\
\lefteqn{ = 2+2 =4= \gamma_{M}(G)   } \hspace*{86mm} \\
\end{eqnarray*}

\vspace{-5mm}

Obviously, it is bigger than $\frac{\beta(G)+2}{3} \ (=
\frac{10}{3}$), and is bigger than
$min\{f_{2}(\delta,g)(\beta(G)+1), \lfloor\frac{\beta(G)}{2}\rfloor
\} \ (=f_{2}(\delta,g)(\beta(G)+1)=3)$.



In the graph $G$ depicted in the right of Fig.10, let $A=\{v_1, v_2,
v_3, v_4, v_5\}$. Then

\vspace{-5mm}

\begin{eqnarray*}
\frac{1}{2}\sum_{i=1}^{m}\big(d(v_{i})-\varepsilon_{i}\big)+
\gamma_{M}\big(G-\{v_1, v_2,  \dots,
 v_{m}\}\big) \\
\lefteqn{ = \frac{1}{2}\big( (3-1)\times5  \big) + \gamma_{M} \big(
G-\{v_1, v_2, v_3, v_4, v_5 \}  \big)   }  \hspace*{86mm} \\
\lefteqn{ = 5+0 =5= \gamma_{M}(G)   } \hspace*{86mm} \\
\end{eqnarray*}

\vspace{-5mm}

Obviously, it is bigger than $\frac{3\beta(G)+4}{7} \ (=
\frac{34}{7})$, and is bigger than
$min\{f_{3}(\delta,g)(\beta(G)+1), \lfloor\frac{\beta(G)}{2}\rfloor
\} \ (=f_{3}(\delta,g)(\beta(G)+1)=\frac{33}{7})$.

\bigskip

\noindent {\bf 4. Independence number and  the maximum genus of
graphs}

\bigskip

Caro\cite{car} and Wei\cite{wei} independently shown that for a
graph $G$ its independence number $$\alpha(G)\geqslant\sum\limits_{
v\in V(G)  }\frac{1}{d_{G}(v)+1}.$$ Later, Alon and Spencer
\cite{alo} gave an elegant probabilistic proof of this bound. But,
up to now, there is little result concerning the relation between
the independence number and the maximum genus of graphs. Let
$N_{G}(v)$ denote all the neighbors of the vertex $v$ in $G$, the
following theorem remedies this deficiency.

\bigskip

{\bf Theorem D }   \ Let $G=(V,E)$ be a connected 3-regular graph
(loops and multi-edges are permitted) with $A=\{x_1, x_2, \dots,
x_{\gamma_{M}(G)}\}$ be a maximum non-separating independent set of
$G$. Then its independence number $$\alpha(G)\geqslant
\gamma_{M}(G)+\alpha(G-N_{A}),$$ where $\alpha(G-N_{A})$ is the
independence number of the subgraph $G-N_{A}$ and $N_{A}$ is the
$closed$ $closure$ of the set $N_{G}\{x_1, x_2, \dots,
x_{\gamma_{M}(G)}\}$, $i.e.$,
$N_{A}=\big(\bigcup_{i=1}^{\gamma_{M}(G)} N_{G}(x_{i}) \big)\cup\{
x_1, x_2, \dots, x_{\gamma_{M}(G)} \}$.

\bigskip

{\bf Proof}   \ From Lemma 3.1 we can get that there exists a
maximum non-separating independent set $A=\{x_1, x_2, \dots,
x_{\gamma_{M}(G)}\}$ which satisfies $G-A$ is connected. Let
$\mathcal {I}$ be an arbitrary independent set of $G-N_{A} $. It is
obvious that every vertex in $A$ is not adjacent to any vertex in
$\mathcal {I}$. So, $A \cup  \mathcal {I} $ is an independent set of
$G$, and the theorem is obtained. $\hspace*{\fill} \Box$

\bigskip

{\bf Remark  }   \  In the graph $G$ depicted in Fig.11, we may
select $A=\{x_1\}$. Then $N_{A}=\{x_1, x_2, x_6\}$, and
$\alpha(G-N_{A})=2$. Noticing $\alpha(G)=3$ and $\gamma_{M}(G)=1$,
we can get that $\alpha(G)= \gamma_{M}(G)+\alpha(G-N_{A})=3>\sum_{
v\in V(G) }\frac{1}{d_{G}(v)+1}=\frac{6}{3+1}=\frac{3}{2}$. So, the
lower bound in Theorem D may be best possible, and may be better
than that of Caro\cite{car} and Wei\cite{wei} in the case of cubic
graphs.

\medskip


\setlength{\unitlength}{1mm}
\begin{center}
\begin{picture}(100,28)

\put(37,7){\begin{picture}(10,10)

\put(0,20){\circle*{2}}  

\put(20,20){\circle*{1.5}}   

\put(35,10){\circle*{2}}   

\put(20,0){\circle*{1.5}}   

\put(0,0){\circle*{2}}   

\put(-15,10){\circle*{1.5}}   

\qbezier(0,20)(-7,13)(-15,10)  

\qbezier(0,20)(-7,19)(-15,10)  

\put(0,20){\line(1,0){20}}   

\qbezier(20,20)(27,13)(35,10)  

\qbezier(20,20)(27,19)(35,10)  

\put(35,10){\line(-3,-2){15}}   

\put(-15,10){\line(3,-2){15}}   

\qbezier(20,0)(10,2.5)(0,0)  

\qbezier(20,0)(10,-2.5)(0,0)  

\begin{footnotesize}

\put(0,23){{$x_1$}}

\put(20,23){{$x_2$}}

\put(-20,10){{$x_6$}}

\put(37,9){{$x_3$}}

\put(-2,-4){{$x_5$}}

\put(21,-4){{$x_4$}}

\end{footnotesize}


\put(5,-9.5){{\bf Fig.11}}

\end{picture}}

\end{picture}
\end{center}

\bigskip

\noindent {\bf 5. Estimating  the number of the maximum genus
embedding of $K_{m}$}

\bigskip

The enumeration of the distinct maximum genus embedding plays an
important role in the study of the genus distribution  problem,
which  may be used to decide whether two given graphs are
isomorphic. But up to now, except \cite{sta} and \cite{ren}, there
is little result concerning the number of the maximum genus
embedding of graphs. In this section, we will provide an algorithm
to enumerate the number of the distinct maximum genus embedding of
the complete graph $K_{m}$, and offer a lower bound which is better
than that of S. Stahl$^{\cite{sta}}$ for $m\leqslant10$.
Furthermore, the enumerative method below can be used to any maximum
genus embedding, other than the method in \cite{sta} which is
restricted to upper embeddable graphs.

A 2-$path$ is called a $\mathcal {V}$-$type$-$edge$, and is denoted
by $\mathcal {V}$. If the $\mathcal {V}$-$type$-$edge$ consists of
the 2-$path$ $v_{i}v_{j}v_{k}$, then this $\mathcal
{V}$-$type$-$edge$ is denoted by  $\mathcal {V}_{j}^{i,k}$ for
simplicity. Let $\psi(G)$ be an embedding of a graph $G$. We say
that a $\mathcal {V}$-$type$-$edge$ are $inserted$ into $\psi(G)$ if
the three endpoints of the $\mathcal {V}$-$type$-$edge$ are inserted
into the corners of the faces in $\psi(G)$, yielding an embedding of
$G+\mathcal {V}$. The following observation can be easily obtained
and is essential in this section.

\medskip

{\bf Observation} \  \ Let $\psi(G)$ be an embedding of a graph $G$.
We can insert a $\mathcal {V}$-$type$-$edge$ $\mathcal {V}$ to
$\psi(G)$ to get an embedding  $\rho(G+\mathcal {V})$ of $G+\mathcal
{V}$ so that the face number of $\rho(G+\mathcal {V})$ is not more
than that of $\psi(G)$.

\medskip

{\bf Lemma 5.1} \  \ Let $\psi(G)$ be a $one$-$face$ embedding of
the graph $G$, $v_{j}$, $v_{i}$ and $v_{k}$ be vertices of $G$. If
the number of the  $face$-$corner$ which containing $v_{j}$, $v_{i}$
and $v_{k}$ are $r_1$, $r_2$ and  $r_3$ respectively, then there are
$r_1 \times r_2 \times r_3$ different ways to add the $\mathcal
{V}$-$type$-$edge$ $\mathcal {V}_{j}^{i, k}$ to $\psi(G)$ to get a
$one$-$face$ embedding of the graph $G+\mathcal {V}_{j}^{i, k}$.

\medskip

{\bf Proof} \ Let the graph depicted in the middle of Fig.12 denote
a $one$-$face$ embedding $\psi(G)$ of the graph $G$. Because the
number of the  $face$-$corner$ which containing $v_{j}$, $v_{i}$ and
$v_{k}$ are $r_1$, $r_2$ and  $r_3$ respectively, we can insert the
$\mathcal {V}$-$type$-$edge$ $\mathcal {V}_{j}^{i, k}$ into
$\psi(G)$ so that there are $r_1$ different ways to put the edges
$v_{j}v_{k}$ and $v_{j}v_{i}$ in the same $face$-$corner$ which
containing the vertex $v_{j}$, $r_2$ different ways to  put the edge
$v_{j}v_{i}$ in a $face$-$corner$ which containing  the vertex
$v_{i}$, and $r_3$ different ways to  put the edge $v_{j}v_{k}$ in a
$face$-$corner$ which containing  the vertex $v_{k}$. For any one of
the $r_1 \times r_2 \times r_3$ different ways to insert the
$\mathcal {V}$-$type$-$edge$ $\mathcal {V}_{j}^{i, k}$ into
$\psi(G)$, we can always get a $one$-$face$ embedding of $G+\mathcal
{V}_{j}^{i, k}$ by one and only one of the two ways which is
depicted by the left and right of Fig.12. So the lemma is obtained.
$\hspace*{\fill} \Box$

\medskip

\setlength{\unitlength}{0.85mm}
\begin{center}
\begin{picture}(100,33)

\put(-28,6){\begin{picture}(10,10)

\qbezier(20,4.7)(10,15)(10,22)

\qbezier(0,10)(20,0)(40,10)

\qbezier(20,5.3)(30,15)(30,21)

\put(20,4.7){\circle*{1.6}}   

\put(20,4.7){\line(-2,-1){10}}

\put(20,4.7){\line(2,-1){10}}

\put(20,-1){\circle*{0.8}}

\put(17,-0.5){\circle*{0.8}}

\put(23,-0.5){\circle*{0.8}}

\put(10,22){\circle*{1.6}}   

\put(10,22){\line(-1,1){6}}

\put(10,22){\line(1,1){6}}

\put(10,22){\line(-3,-2){6}}

\put(10,22){\line(3,-2){6}}

\put(7.5,26.5){\circle*{0.8}}

\put(9.5,27){\circle*{0.8}}

\put(11.7,26.7){\circle*{0.8}}

\put(30,22){\circle*{1.6}}   

\put(30,22){\line(-1,1){6}}

\put(30,22){\line(1,1){6}}

\put(30,22){\line(-3,-2){6}}

\put(30,22){\line(3,-2){6}}

\put(27.3,26.5){\circle*{0.8}}

\put(29.5,27){\circle*{0.8}}

\put(32,26.5){\circle*{0.8}}

\begin{footnotesize}

\put(35,22){{$v_{k}$}}

\put(3,22){{$v_{i}$}}

\put(10,3){{$v_{j}$}}

\end{footnotesize}


\put(72,-11){{\bf Fig.12}}

\end{picture}}

\put(31,6){\begin{picture}(10,10)

\qbezier(0,10)(20,0)(40,10)

\put(20,4.7){\circle*{1.6}}   

\put(20,4.7){\line(-2,-1){10}}

\put(20,4.7){\line(2,-1){10}}

\put(20,-1){\circle*{0.8}}

\put(17,-0.5){\circle*{0.8}}

\put(23,-0.5){\circle*{0.8}}

\put(10,22){\circle*{1.6}}   

\put(10,22){\line(-1,1){6}}

\put(10,22){\line(1,1){6}}

\put(10,22){\line(-3,-2){6}}

\put(10,22){\line(3,-2){6}}

\put(7.5,26.5){\circle*{0.8}}

\put(9.5,27){\circle*{0.8}}

\put(11.7,26.7){\circle*{0.8}}

\put(30,22){\circle*{1.6}}   

\put(30,22){\line(-1,1){6}}

\put(30,22){\line(1,1){6}}

\put(30,22){\line(-3,-2){6}}

\put(30,22){\line(3,-2){6}}

\put(27.3,26.5){\circle*{0.8}}

\put(29.5,27){\circle*{0.8}}

\put(32,26.5){\circle*{0.8}}

\begin{footnotesize}

\put(35,22){{$v_{k}$}}

\put(3,22){{$v_{i}$}}

\put(10,3){{$v_{j}$}}

\end{footnotesize}


\end{picture}}


\put(87,5){\begin{picture}(10,10)

\qbezier(20,5)(25,11)(15,15)

\qbezier(0,10)(20,0)(40,10)

\qbezier(15,15)(9,17)(10,22)

\put(20,4.7){\circle*{1.6}}   

\qbezier(30,21.2)(31,15)(25,13)

\put(20,4.7){\line(-2,-1){10}}

\qbezier(19.2,5.3)(12,10)(25,13)

\put(20,4.7){\line(2,-1){10}}

\put(20,-1){\circle*{0.8}}

\put(17,-0.5){\circle*{0.8}}

\put(23,-0.5){\circle*{0.8}}

\put(10,22){\circle*{1.6}}   

\put(10,22){\line(-1,1){6}}

\put(10,22){\line(1,1){6}}

\put(10,22){\line(-3,-2){6}}

\put(10,22){\line(3,-2){7}}

\put(7.5,26.5){\circle*{0.8}}

\put(9.5,27){\circle*{0.8}}

\put(11.7,26.7){\circle*{0.8}}

\put(30,22){\circle*{1.6}}   

\put(30,22){\line(-1,1){6}}

\put(30,22){\line(1,1){6}}

\put(30,22){\line(-3,-2){6}}

\put(30,22){\line(3,-2){6}}

\put(27.3,26.5){\circle*{0.8}}

\put(29.5,27){\circle*{0.8}}

\put(32,26.5){\circle*{0.8}}

\begin{footnotesize}

\put(35,22){{$v_{k}$}}

\put(3,22){{$v_{i}$}}

\put(10,3){{$v_{j}$}}

\end{footnotesize}

\end{picture}}
\end{picture}
\end{center}

\bigskip

The following algorithm together with Lemma 5.1 provide a maximum
genus embedding of $K_{m}$ and a  lower bound of the number of the
maximum genus embedding of $K_{m}$.

\bigskip

{\bf Algorithm }

\medskip

\textbf{Note:} \ Let $V=\{v_1, v_2, \dots, v_{m}\}$ be the vertex
set of the complete graph $K_{m}$. In the following algorithm,
$\forall \ i \in \{k, a, b\} \subseteq \{1, 2, \dots, m\}$, if
$i\equiv 0 \ (mod \ m)$, then let $i=m$.

\textbf{Step 1.} Embed the tree $v_2v_3\dots v_{m} v_1$ on the
plane.

\textbf{Step 2.}  Let $k=1$, $a=2$, $b=3$.

\textbf{Step 3.} If the $one$-$face$ embedding of the complete graph
$K_{m}$ is obtained, then stop. Otherwise, go to Step 4.

\textbf{Step 4.} If there are only two vertices  $v_{k}$ and $v_{a}$
that are not adjacent,  then connect them to get a $two$-$face$
embedding of the complete graph $K_{m}$ and stop. Otherwise, go to
Step 5.

\textbf{Step 5.} If there is no edge connecting the vertex $v_{k}$
and $v_{a}$ then go to Step 6. Otherwise, go to Step 10.

\textbf{Step 6.} If any pair of $\{ v_{k}, v_{a}, v_{b}\}$ are not
the same, and there is no edge connecting the vertex $v_{k}$ and
$v_{b}$ then add the $\mathcal {V}$-$type$-$edge$ $\mathcal
{V}_{k}^{a, b}$ to the graph to get a $one$-$face$ embedding  and go
to Step 9. Otherwise, let $b\equiv b+1 \ (mod \ m)$ and go to Step
7.

\textbf{Step 7.} If $b\equiv k-1 \ (mod \ m)$ then go to Step 8.
Otherwise, go back to Step 6.

\textbf{Step 8.} Let $c=k$, $k=a$, $a=c$ ($i.e.,$ exchange $k$ and
$a$). Then go back to Step 3.

\textbf{Step 9.} Let $b\equiv a+3 \ (mod \ m)$, $a\equiv a+2 \ (mod
\ m)$, and go to Step 11.

\textbf{Step 10.} Let $a\equiv a+1 \ (mod \ m)$, and go to Step 11.

\textbf{Step 11.} If $a\equiv k-1 \ (mod \ m)$, then go to Step 12.
Otherwise, go back to Step 3.

\textbf{Step 12.} Let $k=1$, and go  to Step 13.

\textbf{Step 13.} If $d_{G}(v_{k})< m-1$, then let $a\equiv k+2 \
(mod \ m)$, $b\equiv k+3 \ (mod \ m)$, and go back to Step 3.
Otherwise, go to Step 14.

\textbf{Step 14.} If $d_{G}(v_{k})= m-1$, then let $k=k+1$, and go
back to Step 13.

\bigskip

Using the above algorithm, we can get the maximum genus embedding of
$K_{m}$ except that $m=1+8i$ or $m=6+8i$ ($i=0,1,3, \dots$).
Furthermore, for $m\leqslant10$, our result is much better than that
of Stahl$^{\cite{sta}}$. For simplicity, we give some symbols which
are used below. Let $E$ be a $one$-$face$ embedding  of a graph.
Then the symbol ($\mathcal {V}_{j}^{i,k}: r_1\times r_2\times r_3$)
means that there are $r_1\times r_2\times r_3$ different ways to add
the $\mathcal {V}$-$type$-$edge$ $\mathcal {V}_{j}^{i,k}$ to $E$ to
get a $one$-$face$ embedding  of $E+\mathcal {V}_{j}^{i,k}$, and the
symbol ($e_{j}^{j,k}: r_1\times r_2$) means that there are
$r_1\times r_2$ different ways to add the edge $v_{j}v_{k}$ to $E$
to get a $two$-$face$ embedding of $E+v_{j}v_{k}$.

\medskip

{\bf Result 1 } \  \ The number of the maximum genus embedding of
the complete graph $K_8$ is at least $2^{26}\times3^{11}\times5^5$.

\medskip

{\bf Proof} \ \ Let $V=\{v_1, v_2, \dots, v_{8}\}$ be the vertex set
of the complete graph $K_{8}$. There is only one way to embed the
tree $T=v_2v_3\dots v_8v_1$ on the plane, which is a $one$-$face$
embedding, and is denoted by $\mathcal {E}_1$. In $\mathcal {E}_1$,
the number of the $face$-$corner$ which containing the vertex $v_1$,
$v_2$, $v_3$ is 1, 1 and 2 respectively. So, according to Lemma 5.1,
there are 2 different ways to add the $\mathcal {V}$-$type$-$edge$
$\mathcal {V}_{1}^{2,3}$ to $\mathcal {E}_1$ to get a $one$-$face$
embedding  of $T+\mathcal {V}_{1}^{2,3}$. Let $\mathcal {E}_2$ be
any one of the $one$-$face$ embedding of $T+\mathcal {V}_{1}^{2,3}$.
In $\mathcal {E}_2$, the number of the $face$-$corner$ which
containing the vertex $v_1$, $v_4$, $v_5$ is 3, 2 and 2
respectively. So, according to Lemma 5.1, there are $3\times2\times2
\ (=12)$ different ways to add the $\mathcal {V}$-$type$-$edge$
$\mathcal {V}_{1}^{4,5}$ to $\mathcal {E}_2$ to get a $one$-$face$
embedding  of $T+\mathcal {V}_{1}^{2,3}+\mathcal {V}_{1}^{4,5}$.
Similarly, we can get that for each of the $one$-$face$ embedding of
$T+\mathcal {V}_{1}^{2,3}+\mathcal {V}_{1}^{4,5}$, there are
$5\times2\times2$ different ways to add the $\mathcal
{V}$-$type$-$edge$ $\mathcal {V}_{1}^{6,7}$ to $T+\mathcal
{V}_{1}^{2,3}+\mathcal {V}_{1}^{4,5}$ to get a $one$-$face$
embedding  of $T+\mathcal {V}_{1}^{2,3}+\mathcal
{V}_{1}^{4,5}+\mathcal {V}_{1}^{6,7}$.

Similarly, we can add $\mathcal {V}$-$type$-$edges$, one by one in
the following order, to $T+\mathcal {V}_{1}^{2,3}+\mathcal
{V}_{1}^{4,5}+\mathcal {V}_{1}^{6,7}$ to get a $two$-$face$
embedding  of $K_8$ eventually.

($\mathcal {V}_{2}^{4,5}: 2\times 3\times 3$), ($\mathcal
{V}_{2}^{6,7}: 4\times 3\times 3$), ($\mathcal {V}_{8}^{2,3}:
2\times 6\times 3$), ($\mathcal {V}_{8}^{4,5}: 4\times 4\times 4$),
($\mathcal {V}_{6}^{8,3}: 4\times 6\times 4$), ($\mathcal
{V}_{4}^{6,7}: 5\times 6\times 4$), ($\mathcal {V}_{3}^{5,7}:
5\times 5\times 5$), ($e_{5}^{5,7}: 6\times 6$).

So, the number of the distinct maximum genus embedding of $K_8$ is
at least

\vspace{-6mm}

\begin{eqnarray*}
2\times(3\times2\times2)\times(5\times2\times2)\times(2\times3\times3)\times(4\times3\times3)\times(2\times6\times3)  \\
\lefteqn{
\times(4\times4\times4)\times(4\times6\times4)\times(5\times6\times4)\times(5\times5\times5)\times(6\times6)
 }  \hspace*{127mm} \\
\lefteqn{ =2^{26}\times3^{11}\times5^5
 }  \hspace*{132mm} \\
\end{eqnarray*}

\vspace{-6mm}

{\bf Result 2 } \  \ The number of the distinct maximum genus
embedding of the complete graph $K_{10}$ is at least
$2^{52}\times3^{15}\times5^{7}\times7^6$, which is obtained from the
unique $one$-$face$ embedding  of the tree $T=v_2v_3\dots v_{10}v_1$
by successively adding the following $\mathcal {V}$-$type$-$edges$:
($\mathcal {V}_{1}^{2,3}: 1\times 1\times 2$), ($\mathcal
{V}_{1}^{4,5}: 3\times 2\times 2$), ($\mathcal {V}_{1}^{6,7}:
5\times 2\times 2$), ($\mathcal {V}_{1}^{8,9}: 7\times 2\times 2$),
($\mathcal {V}_{2}^{4,5}: 2\times 3\times 3$), ($\mathcal
{V}_{2}^{6,7}: 4\times 3\times 3$), ($\mathcal {V}_{2}^{8,9}:
6\times 3\times 3$), ($\mathcal {V}_{10}^{2,3}: 2\times 8\times 3$),
($\mathcal {V}_{10}^{4,5}: 4\times 4\times 4$), ($\mathcal
{V}_{10}^{6,7}: 6\times 4\times 4$), ($\mathcal {V}_{8}^{10,3}:
4\times 8\times 4$), ($\mathcal {V}_{8}^{4,5}: 6\times 5\times 5$),
($\mathcal {V}_{6}^{8,9}: 5\times 8\times 4$), ($\mathcal
{V}_{6}^{3,4}: 7\times 5\times 6$), ($\mathcal {V}_{3}^{5,7}:
6\times 6\times 5$), ($\mathcal {V}_{9}^{3,4}: 5\times 8\times 7$),
($\mathcal {V}_{9}^{5,7}: 7\times 7\times 6$), ($\mathcal
{V}_{7}^{4,5}: 7\times 8\times 8$).

\medskip

{\bf Result 3 } \  \ The number of the distinct maximum genus
embedding of the complete graph $K_7$ is at least 49766400000, which
is obtained from the unique $one$-$face$ embedding  of the tree
$T=v_2v_3\dots v_7v_1$ by successively adding the following
$\mathcal {V}$-$type$-$edges$: ($\mathcal {V}_{1}^{2,3}: 1\times
1\times 2$), ($\mathcal {V}_{1}^{4,5}: 3\times 2\times 2$),
($\mathcal {V}_{6}^{1,2}: 2\times 5\times 2$), ($\mathcal
{V}_{6}^{3,4}: 4\times 3\times 3$), ($\mathcal {V}_{2}^{4,5}:
3\times 4\times 3$), ($\mathcal {V}_{7}^{2,3}: 2\times 5\times 4$),
($\mathcal {V}_{7}^{4,5}: 4\times 5\times 4$), ($e_{3}^{3,5}:
5\times 5$).

\medskip

{\bf Result 4 } \  \ The number of the distinct maximum genus
embedding of the complete graph $K_5$ is at least 432, which is
obtained from the unique $one$-$face$ embedding  of the tree
$T=v_2v_3v_4v_5v_1$ by successively adding the following $\mathcal
{V}$-$type$-$edges$: ($\mathcal {V}_{1}^{2,3}: 1\times 1\times 2$),
($\mathcal {V}_{4}^{1,2}: 2\times 3\times 2$), ($\mathcal
{V}_{5}^{2,3}: 2\times 3\times 3$).

\medskip

The algorithm doesn't  work for $K_6$ and $K_9$. But the maximum
genus embedding of $K_6$ and $K_9$ can be obtained by the following
manners.

\medskip

{\bf Result 5 } \  \ The number of the distinct maximum genus
embedding of the complete graph $K_6$ is at least 663552, which is
obtained from the unique $one$-$face$ embedding  of the tree
$T=v_2v_3\dots v_6v_1$ by successively adding the following
$\mathcal {V}$-$type$-$edges$: ($\mathcal {V}_{1}^{2,3}: 1\times
1\times 2$), ($\mathcal {V}_{1}^{4,5}: 3\times 2\times 2$),
($\mathcal {V}_{2}^{4,5}: 2\times 3\times 3$), ($\mathcal
{V}_{6}^{2,4}: 2\times 4\times 4$), ($\mathcal {V}_{3}^{5,6}:
3\times 4\times 4$).

\medskip

{\bf Result 6 } \  \ The number of the distinct maximum genus
embedding of the complete graph $K_9$ is at least
$2^{27}\times3^{12}\times5^{7}\times7^6$, which is obtained from the
unique $one$-$face$ embedding  of the tree $T=v_2v_3\dots v_9v_1$ by
successively adding the following $\mathcal {V}$-$type$-$edges$:
($\mathcal {V}_{1}^{2,3}: 1\times 1\times 2$), ($\mathcal
{V}_{1}^{4,5}: 3\times 2\times 2$), ($\mathcal {V}_{1}^{6,7}:
5\times 2\times 2$), ($\mathcal {V}_{8}^{1,2}: 2\times 7\times 2$),
($\mathcal {V}_{8}^{3,4}: 4\times 3\times 3$),   ($\mathcal
{V}_{8}^{5,6}: 6\times 3\times 3$), ($\mathcal {V}_{2}^{4,5}:
3\times 4\times 4$),  ($\mathcal {V}_{2}^{6,7}: 5\times 4\times 3$),
($\mathcal {V}_{9}^{2,3}: 2\times 7\times 4$), ($\mathcal
{V}_{9}^{4,5}: 4\times 5\times 5$), ($\mathcal {V}_{9}^{6,7}:
6\times 5\times 4$), ($\mathcal {V}_{3}^{5,6}: 5\times 6\times 6$),
($\mathcal {V}_{7}^{3,5}: 5\times 7\times 7$), ($\mathcal
{V}_{4}^{6,7}: 6\times 7\times 7$).

\medskip

{\bf Remark}   \ \ Saul Stahl$^{\cite{sta}}$ obtained that the
complete graph $K_{m}$ on $m$ vertices has at least
$[(m-6)!]^4[(m-3)!]^{m-4}$ maximum genus embeddings, and for
$m\equiv 0,3 \ (mod \ 4)$ $K_{m}$ has at least
$(\frac{m-2}{m-1})^2[(m-3)!]^{m}$ maximum genus embeddings. It is
obvious that our results for $m\leqslant10$ is much better than that
of Stahl.


$\bf{Acknowledgements}$
 \ \ The authors  thank the referees for their careful reading of the
paper, and for their valuable comments.

\medskip


\end{document}